\documentclass{elsart}

\usepackage{amsfonts,amssymb}
\newtheorem{corollary}{Corollary}
\newtheorem{theorem}{Theorem}
\newtheorem{lemma}{Lemma}
\newtheorem{proposition}{Proposition}

\newcommand{\pp}{\noindent {\em Proof. }}
\newcommand{\bee}[1]{\begin{equation}\label{#1}}
\newcommand{\beq}[1]{\begin{eqnarray}\label{#1}}
\newcommand{\ene}{\end{equation}}
\newcommand{\eqe}{\end{eqnarray}}
\newcommand{\ld}{\ldots}

\newcommand{\cd}{\cdots}
\newcommand{\tr}[1]{ ^t\!#1}
\newcommand{\rsl}{\mathrm{sl}}

\newcommand{\wh}[1]{\widehat{#1}}
\newcommand{\wg}{\widehat{G}}

\newcommand{\vp}{\varphi}
\newcommand{\bx}{\hfill$\Box$}
\newcommand{\cj}[2]{#1^{-1}\!#2#1}
\newcommand{\iv}[1]{#1^{-1}\!}
\newcommand{\ve}{\varepsilon}

\newcommand{\be}{\bar{e}}
\newcommand{\bg}{\bar{g}}
\newcommand{\bG}{\bar{G}}

\newcommand{\epf}{\hfill$\Box$}
\newcommand{\su}[1]{\mathrm{Supp}\,#1}
\newcommand{\Sp}[1]{\mathrm{Span}\,\{#1\}}
\begin{document}

\begin{frontmatter}



\title{Group Gradings on Simple Lie Algebras of Type ``A''}


\author[yb]{Y. A. Bahturin}
\address{Department of Mathematics and Statistics\\Memorial
University of Newfoundland\\ St. John's, NL, A1C5S7,
Canada and \\ Department of
Algebra, Faculty of Mathematics and Mechanics\\Moscow
State University\\Moscow, 119992, Russia}
\ead{yuri@math.mun.ca}\thanks[yb]{Work is partially supported by NSERC grant \#
227060-04 and URP grant, Memorial University of Newfoundland}

\author[mz]{M. V. Zaicev}
\address{Department of Algebra,
Faculty of Mathematics and Mechanics\\Moscow State
University\\Moscow, 119992, Russia}\ead{zaicev@mech.math.msu.su}
\thanks[mz]{Work is partially supported by
RFBR, grant 02-01-00219, and SSC-1910.2003.1}

\begin{abstract}
In this paper we describe all group gradings by a finite abelian group
$G$ of any Lie algebra $L$ of the type ``A'' over
algebraically closed field $F$ of characteristic zero.
\end{abstract}

\begin{keyword}
Graded algebra \sep simple Lie algebra \sep matrix algebra\sep involution
\end{keyword}
\end{frontmatter}



\section{Introduction}\label{s0}

In this paper we describe all gradings by finite abelian groups on the matrix Lie algebras $\rsl(n)$ over an algebraically closed field of characteristic zero. We introduce two types of gradings of $\rsl(n)$, type I, induced from the gradings of the full matrix algebras $M_n$ , described in \cite{BSZ}, and type II, obtained by a simple procedure from so called involution gradings of the full matrix algebras $M_n$, described in \cite{BShZ}, and show that they exhaust all possible gradings of $\rsl(n)$. A survey of known results about the gradings of simple Lie algebras can be found in \cite[Section 3.3]{OV}, including the results of V. Kac \cite{VK} classifying cyclic gradings of all simple Lie algebras. Many examples of grading on simple Lie algebras can be found in the papers of J. Patera and coauthors (see, for example, \cite{HPPT,JP1,JP2}). 

\section{Some notation and simple facts}\label{s1}
Let $F$ be an arbitrary field, $A$ a not necessarily associative algebra over an $F$ and $G$ a group. We say
that  $A$ is a $G$-graded algebra, if there is a vector space sum
decomposition
\bee{e0001}
A=\bigoplus_{g\in G} A_g,
\ene
such that
\bee{e0002}
A_gA_h\subset A_{gh}\mbox{ for all }g,h\in G.
\ene

A subspace $V\subset A$ is called {\em graded} (or {\it homogeneous})
if $V=\oplus_{g\in G} (V\cap A_g)$. An element $a\in R$ is called
{\it homogeneous of degree} $g$ if $a\in A_g$. We also write $\deg a=g$. The {\it support} of the
$G$-grading is a
subset
$$
\mathrm{Supp}\: A=\{g\in G|A_g\ne 0\}.
$$

If $N$ is a normal subgroup of $G$ then $A$ naturally acquires a $G/N$-grading if one sets

\bee{efg}
A_x=\bigoplus_{g\in x} A_g\mbox{ for any }x\in G/N.
\ene

If $A$ is an associative algebra then (\ref{e0001}) is called a \textit{Lie grading} if instead of (\ref{e0002}) one has 
\bee{e0003}
[A_g,A_h]\subset A_{gh}\mbox{ for all }g,h\in G.
\ene

Suppose now that $F$ is of characteristic different from $2$. If $A$ is an associative algebra with involution $\ast$ and, in addition to (\ref{e0002}),  one has
\bee{e0004}
(A_g)^\ast\subset A_g\mbox{ for all }g\in G.
\ene
then we say that (\ref{e0001}) is an \textit{involution preserving} grading or simply an \textit{involution} grading. In this case, given a graded subspace $B\subset A$ we set
\bee{esym}
B^{(+)}=\{ b\in B\,|\,b^*=b\},\mbox{ the set of symmetric elements of }B
\ene
and
\bee{essym}
B^{(-)}=\{ b\in B\,|\,b^*=-b\},\mbox{ the set of skew-symmetric elements of }B.
\ene
If $B$ is an associative subalgebra of $A$ then $B^{(-)}$ is a Lie subalgebra of $A$, that is, with respect to $[x,y]=xy-yx$ while $B^{(+)}$ is a Jordan subalgebra of $A$, that is, with respect to $x\circ y=xy+yx$. We always have $B=B^{(-)}\oplus B^{(+)}$.

In this paper we study group gradings on \textit{simple Lie algebras}. The following simple remark formally shows why in this case we may restrict ourselves to the case of abelian groups. We give a proof for completeness.

\begin{lemma}\label{l0001}
Let $L=\oplus_{g\in G}L_g$ be a simple Lie algebra over an arbitrary field $F$, graded by a (possibly, noncommutative) group $G$.
Then the support $\mathrm{Supp}\,L$ generates in $G$ an abelian subgroup.
\end{lemma}

\pp First we note that $gh=hg$ for $g,h\in G$ as soon as
$[L_g,L_h]\ne 0$ since $[L_g,L_h]=[L_h,L_g]\subset L_{gh}\cap
L_{hg}$. We are going to generalize this property to the case of more than two
factors. Namely, we will prove that inequality
$[L_{g_1},\ldots,L_{g_k}]\ne 0$ implies that $g_1,\ldots,g_k$
pairwise commute.

Suppose $k\ge 3$ and $[x_1,\ldots,x_k]\ne 0$ for some $x_1\in
L_{g_1},\ldots,x_k\in L_{g_k}$. Then $y=[x_1,\ldots,x_{k-1}]\ne 0$
and by induction all $g_1,\ldots,g_{k-1}$ commute. Also $g_k$
commutes with the product $g_1\cdots g_{k-1}$. Clearly, at least
one of two products $a=[x_1,\ldots,x_{k-2},x_k]$,
$b=[[x_1,\ldots,x_{k-2}],[x_{k-1},x_k]]$ is non-zero. If $a\ne 0$
then $g_k$ commutes with all $g_1,\ldots,g_{k-2}$ and with
$g_1\cdots g_{k-1}$. Hence $g_k$ commutes with $g_{k-1}$.
Similarly, if $b\ne 0$ then $g_{k-1}g_{k}=g_{k}g_{k-1}$ commutes
with all $g_1,\ldots,g_{k-2}$. Hence $g_k$ commutes with all
$g_1,\ldots,g_{k-1}$.

Now we consider arbitrary $g,h\in \mathrm{Supp}\,L$. Since $L$ is simple and
$L_g\ne 0$, there exist $g_1,\ldots,g_k\in \mathrm{Supp}\,L$ such that
$[L_g,L_{g_1},\ldots, L_{g_k}]\ne 0$ and $gg_1\cdots g_k=h$. Then
all $g,g_1,\ldots,g_k$ commute and hence $gh=hg$.\epf

\section{Two types of Lie gradings on associative algebras}\label{stwotypes}

If $R$ is an associative algebra over an arbitrary field $F$, graded by an {\it abelian} group $G$, $L$ a Lie subalgebra in $R$ with respect to the bracket operation $[x,y]=xy-yx$ and $L$ is a graded subspace of $R$ then $L$ becomes a $G$-graded algebra with $L_g=L\cap R_g$. The inclusion $[L_g,L_h]\subset L_{gh}$ easily follows. Thus some gradings of a Lie algebra can be induced from the gradings of an associative algebra. In certain cases all gradings of important Lie algebras can be obtained in this way.

For example, this is the case when $R=M_n$, the matrix algebra of order $n$ over an algebraically closed field $F$ of characteristic 0 and $L$ is a Lie subalgebra of all matrices which are skew-symmetric under a symplectic involution (simple Lie algebra of the type $C_k$, $n=2k$, $k\geq 3$) or a Lie subalgebra  of all matrices which are skew-symmetric under a transpose involution (simple Lie algebra of the type $B_k$, $n=2k+1$, $k\geq 2$ or simple Lie algebra of the type $D_k$, $n=2k$, $k>4$). This made it possible to give a complete description of all abelian gradings on simple Lie algebras of the types just mention, in an earlier paper \cite{BShZ}.

Below we briefly recall the results of \cite{BSZ}, where the full description of abelian group gradings on the full matrix algebra has been given.

\subsection{Abelian Gradings on Matrix Algebras}\label{sAGMA}

A grading $R=\oplus_{g\in G}R_g$ on
the matrix algebra $R=M_n(F)$ is called {\it elementary} if there exists an $n$-tuple
$(g_1,\ldots,g_n)\in G^n$ such that the matrix units $E_{ij}, 1\le
i,j\le n$ are homogeneous and $ E_{ij}\in R_g\iff g=g_i^{-1}g_j. $

A grading is called {\it fine} if $\dim R_g=1$ for any $g\in
\mathrm{Supp}\: R$. A particular case of fine gradings is the so-called
$\varepsilon$-grading where $\varepsilon$ is $n^\mathrm{th}$ primitive root
of $1$. Let $G=\langle a\rangle_n\times\langle b\rangle_n$ be the direct product of two cyclic groups of order
$n$ and

\bee{e00}
X_a = \left(\begin{array}{cccc} \varepsilon^{n-1} & 0 & ... & 0
\\ 0 & \varepsilon^{n-2} & ... & 0
\\  \cdots& \cdots & \cdots& \cdots
\\ 0 & 0 & ... & 1
\end{array}\right)~,~~
X_b = \left(\begin{array}{cccc}
   0 & 1 &  ... & 0
\\ \cdots & \cdots   &    \cdots &  \cdots
\\ 0 & 0 &  ... & 1
\\ 1 & 0 &  ... & 0
\end{array}\right)~.
\ene
Then
\begin{equation}\label{a1}
X_aX_bX_a^{-1}=\varepsilon X_b~,~~X_a^n=X_b^n=I
\end{equation}
and all $X_a^iX_b^j, 1\le i,j\le n$, are linearly independent.
Clearly, the elements $X_a^iX_b^j,\, i,j=1,\ldots, n$, form a basis
of $R$ and all the products of these basis elements are uniquely
defined by (\ref{a1}).

        Now for any $g\in G, g=a^ib^j$, we set $X_g=X_a^iX_b^j$ and denote by $R_g$
a one-dimensional subspace
\bee{aa1}
R_g=\langle X_a^iX_b^j\rangle.
\ene
Then from (\ref{a1}) it follows that $R=\oplus_{g\in G}R_g$ is a
$G$-grading on $M_n(F)$ which is called an $\varepsilon$-grading.

Now let $R=M_n(F)$ be the full matrix algebra over $F$ graded by an abelian
group $G$. The folowing result has been proved in \cite[Section 4, Theorems 5, 6]{BSZ} and \cite[Subsection 2.2, Theorem 6]{surgrad}.

\begin{theorem}\label{tbsz} Let $F$ be an algebraically closed field of characteristic zero. Then as a $G$-graded algebra $R$ is isomorphic to the
tensor product
$$
R^{(0)}\otimes R^{(1)}\otimes\cd\otimes R^{(k)}
$$
where $R^{(0)}=M_{n_0}(F)$ has an elementary $G$-grading, $\mathrm{Supp}\:
R^{(0)}=S$ is a
finite subset of $\,$ $G$, $R^{(i)}=M_{n_i}(F)$ has the $\varepsilon_i$ grading,
$\varepsilon_i$ being a primitive $n_i^\mathrm{th}$ root of $1$,
$\mathrm{Supp}\; R^{(i)}=H_i\cong \mathbb{Z}_{n_i}\times
\mathbb{Z}_{n_i}, i=1,\ld,k$. Also $H=H_1\cd H_k\cong H_1\times\cd\times
H_k$ and $S\cap H=\{e\}$ in $G$. 
\end{theorem}

Since $L=\rsl(n)=[R,R]$, $L$ is a graded subspace of $R=M_n$, and the following is true.

\begin{corollary}\label{cinn}\textsc{(Type I Gradings)} Given any grading $R=\bigoplus_{g\in G} R_g$ by a finite abelian group $G$, setting $L_g=R_g\cap L$ makes $L$ into a $G$-graded Lie algebra.
\end{corollary}

If $R=M_n$ and $L=\rsl(n)\subset M_n$ then not all gradings of $L$ are induced from $R$. For example, if $n>2$, $G=\mathbb{Z}_2$, $L_0$ is the set of all skew-symmetric matrices under the ordinary transpose involution $X\rightarrow    \;^t\! X$, $L_1$ the set of all symmetric matrices of trace zero. This is not induced from any $\mathbb{Z}_2$-grading of $M_n$ since, according to the general theory of \cite[Section 4, Theorems 5, 6]{BSZ} any such grading is elementary, that is, there are two natural numbers $k,l$ such that $\dim R_0 = k^2+l^2$ and $\dim R_1=2kl$. In this case $\dim L_0=k^2+l^2-1=\frac{n(n-1)}{2}$, and $\dim L_1=2kl=\frac{n(n+1)}{2}-1$. Since $k^2+l^2\geq 2kl$ we easily derive $n\leq 2$.

Quite a general result, as we will see in the future, dealing with the gradings of an associative algebra $R$ in the presence of involutions, is the following.

\begin{theorem}\label{tomega}Let $R$ be an associative algebra over a field $F$ of characteristic different from $2$, $G$ be a finite abelian group, $h$ an element of order 2 in $G$. If $\ast$ is an involution on $R$ and
$$
R=\bigoplus_{g\in G}\widetilde{R}_g
$$
is an involution $G$-grading then a  Lie grading by $G$ on $R$ can be given by
$$
R_g=\widetilde{R}_g^{(-)}\oplus \widetilde{R}_{gh}^{(+)}.
$$
Here, as in Section \ref{s1}, $\widetilde{R}_g^{(\pm)}$ is the set of symmetric (skew-symmetric) elements in $\widetilde{R}_g$ with respect to $\ast$.\epf
\end{theorem}

The proof of this theorem is a direct verification of the relations $[R_g,R_{g^\prime}]\subset R_{gg^\prime}$, for all $g,g^\prime\in G$. If we consider the restriction of the grading of Theorem \ref{tomega} on $R=M_n$ to the matrices of trace zero then we obtain the following result.

\begin{corollary}\label{cout}\textsc{(Type II Gradings)} Let $R=M_n$ and suppose $R$ satisfies the hypotheses of Theorem \ref{tomega}. If also $\mathrm{char}\,F \nmid n$ then a grading of the Lie algebra $L=\rsl(n)$ will be obtained if we set
$$
L_g=\left\{\begin{array}{ll}\widetilde{R}_g^{(-)}\oplus \widetilde{R}_{gh}^{(+)}&\mbox{ if }g\neq h\\
\widetilde{R}_h^{(-)}\oplus (\widetilde{R}_{e}^{(+)}\cap L)&\mbox{ otherwise}.\end{array}\right.
$$
\end{corollary}

Note that we need the restriction $\mathrm{char}\,F \nmid n$ to make sure that $R = FI \oplus [R,R]$.

All involution gradings of $R=M_n$ have been described in \cite[Sections 6, 7, 8]{BShZ}.

\subsection{Abelian Involution Gradings on Matrix Algebras}\label{sAIGMA} Throughout this subsection we assume the base field $F$ being algebraically closed of characteristic zero.

We first recall that any involution $\ast$ of $R=M_n$ can always be written as

\bee{einv}
X^\ast= \Phi^{-1}( \tr{X})\Phi
\ene
where $\Phi$ is a nondegenerate matrix which is either symmetric or skew-sym\-metric and $X\mapsto\, \tr{X}$ is the ordinary transpose map. In the case where $\Phi$ is symmetric we call $\ast$ a transpose involution. If $\Phi$ is skew-symmetric $\ast$ is called a symplectic involution. Before we formulate the theorem describing involution gradings on $M_n$ we need three (slightly modified) lemmas from \cite{BShZ}. The first one deals with certain fine involution gradings while the last two with elementary involution gradings. If $R$ has an involution $\ast$ then by $R^{(\pm)}$ we denote the space of symmetric (respectively skew-symmetric) matrices in $R$ under $\ast$.

\begin{lemma}\label{L6}
Let $R=M_2(F)$ be a $2\times 2$ matrix algebra endowed with an involution $*:R\to R$ corresponding to a symmetric or skew-symmetric
non-degenerate bilinear form with the matrix $\Phi$. Then the $(-1)$-grading
of $M_2$ is an involution grading if and only if one of the following
holds:
\begin{itemize}
\item[(1)] $\Phi$ is skew-symmetric,
$$
\Phi=\pmatrix{0 & 1 \cr -1 & 0 \cr},
\quad
R^{(-)}=\left\{ \pmatrix{ a & b \cr c & -a \cr}\right\},
\quad
R^{(+)}=\left\{ \pmatrix{ a & 0 \cr 0 & a \cr}\right\};
$$
and
$$
\pmatrix{x & y \cr z & t \cr}^* = \pmatrix{t & -y \cr -z & x \cr};
$$

\item[(2)] $\Phi$ is symmetric,
$$
\Phi=\pmatrix{0 &  1 \cr  1 & 0 \cr},
\quad
R^{(-)}=\left\{ \pmatrix{ a & 0 \cr 0 & -a \cr}\right \},
\quad
R^{(+)}=\left\{ \pmatrix{ a & b \cr c & a \cr}\right \};
$$
and
$$
\pmatrix{x & y \cr z & t \cr}^* = \pmatrix{t & y \cr z & x \cr};
$$

\item[(3)] $\Phi$ is symmetric,
$$
\Phi=\pmatrix{1 &  0 \cr  0 & 1 \cr},
\quad
R^{(-)}=\left\{ \pmatrix{ 0 & b \cr -b & 0 \cr}\right\},
\quad
R^{(+)}=\left\{ \pmatrix{ a & b \cr b & c \cr}\right\};
$$
and
$$
\pmatrix{x & y \cr z & t \cr}^* = \pmatrix{x & z \cr y & t \cr};
$$

\item[(4)] $\Phi$ is symmetric,
$$
\Phi=\pmatrix{1 &  0 \cr  0 & -1 \cr},
\quad
R^{(-)}=\left\{ \pmatrix{ 0 & b \cr b & 0 \cr}\right\},
\quad
R^{(+)}=\left\{ \pmatrix{ a & b \cr -b & c \cr}\right\};
$$
and
$$
\pmatrix{x & y \cr z & t \cr}^* = \pmatrix{x & -z \cr -y & t \cr}.
$$
\end{itemize}
\end{lemma}
The next lemma handles the case of an elementary grading compatible with an involution defined by a symmetric non-degenerate bilinear form.

\begin{lemma}\label{L8'}
Let $R=M_n(F)$, $n$ a natural number, be a matrix algebra with involution *
defined by a symmetric non-degenerate bilinear form. Let $G$ be an
abelian group and let $R$ be equipped with an elementary
involution $G$-grading defined by an $n$-tuple $(g_1,\ld,g_n)$. Then
$g_1^2=\ld=g_m^2=g_{m+1}g_{m+l+1}=\ld=g_{m+l}g_{m+2l}$ for some
$0\le l\le\frac{n}{2}$ and $m+2l=n$. The involution * acts as
$X^*=(\Phi^{-1})\:^tX\Phi$ where
$$
 \Phi=\pmatrix{   I_m   &   0     &  0   \cr
                   0    &   0     &  I_l \cr
                   0    &   I_l   &  0   \cr},
$$
where $I_s$ is the $s\times s$ identity matrix. Moreover, $R^{(-)}$
consists of all matrices of the type
\bee{BD1}
       \pmatrix{  P    &   S     &   T     \cr
                -^t T   &   A     &   B     \cr
                -^t S   &   C     & -^t A    \cr},
\ene
where $^t P=-P,~ ^t B=-B,~ ^t C=-C$ and
$$
P\in M_m(F),~ A,B,C,D\in M_l(F),~ S,T\in M_{m\times l}(F)
$$
while $R^{(+)}$ consists of all matrices of the type
\bee{BD2}
       \pmatrix{  P    &   S     &   T    \cr
                 ^t T   &   A     &   B    \cr
                 ^t S   &   C     &  ^t A   \cr},
\ene
where $^t P= P,~ ^t B= B,~ ^t C= C$ and
$$
P\in M_m(F),~ A,B,C,D\in M_l(F),~ S,T\in M_{m\times l}(F).
$$
\end{lemma}

The last lemma deals with the case of an elementary grading compatible with an involution defined by a skew-symmetric non-degenerate bilinear form.

\begin{lemma}\label{L8}
Let $R=M_n(F)$, $n=2k$, be the matrix algebra with involution $\ast$
defined by a skew-symmetric non-degenerate bilinear form. Let $G$ be
an abelian group and let $R$ be equipped with an elementary
involution $G$-grading defined by an $n$-tuple $(g_1,\ld,g_n)$. Then $g_{1}g_{k+1}=\ld=g_{k}g_{2k}$,
the involution $\ast$ acts as $X^*=(\Phi^{-1})\:^tX\Phi$ where
$$
 \Phi=\pmatrix{  0   &   I     \cr
                 -I    & 0    \cr},
$$
$I$ is  the $k\times k$ identity matrix, $R^{(-)}$ consists of all
matrices of the type
\bee{C1}
       \pmatrix{  A    &   B     \cr
                  C    & -^t A    \cr},~
A,B,C,\in M_k(F),~ ^t B=B,^t C=C
\ene
while $R^{(+)}$ consists of all matrices of the type
\bee{C2}
       \pmatrix{  A    &   B     \cr
                  C    &  ^t A    \cr},~
A,B,C,\in M_k(F),~ ^t B=-B,^t C=-C.
\ene
\end{lemma}

\bigskip

We can now explicitly describe gradings on a matrix algebra with
involution. Recall that in the case of an algebraically closed field any involution $\ast$ on a matrix algebra is,  up
to isomorphism, either the transpose
or the symplectic involution. 

\begin{theorem}\label{inv}
Let $R=M_n(F)=\oplus_{g\in G}R_g$ be a matrix algebra over an
algebraically closed field of characteristic zero graded by the
group $G$ and $\mathrm{Supp}\: R$ generates $G$. Suppose that
$*: R\rightarrow R$ is a graded involution. Then $G$ is abelian,
$n=2^km$ and $R$ as a $G$-graded algebra with involution is
isomorphic to the tensor product $R^{(0)}\otimes
R^{(1)}\otimes\cdots\otimes R^{(k)}$ where
\begin{itemize}
\item[(1)]
$R^{(0)},\ldots, R^{(k)}$ are graded subalgebras stable under the involution $\ast$;
\item[(2)]
$R^{(0)}=M_m(F)$ is as in Lemma \ref{L8} if $\ast$ is symplectic on
$R^{(0)}$ or as in Lemma \ref{L8'} if $\ast$ is transpose on $R_0$;
\item[(3)]
$R^{(1)}\otimes\cdots\otimes R^{(k)}$ is a $T=T_1\times\ldots\times
T_k$-graded algebra and any $R^{(i)},1\le i\le k$, is $T_i\cong
{\mathbb Z}_2\times{\mathbb Z}_2$-graded algebra as in Lemma
\ref{L6}.
\item[(4)] A graded basis of $R$ is formed by the elements $Y\otimes X_{t_1}\otimes\cdots\otimes X_{t_k}$, where $Y$ is an element of a graded basis of $R^{(0)}$ and the elements $X_{t_i}$ are of the type (\ref{e00}), with $n=2$, $t_i\in T_i$. The involution on these elements is given canonically by 
$$(Y\otimes X_{t_1}\otimes\cdots\otimes X_{t_k})^\ast=Y^\ast\otimes X_{t_1}^\ast\otimes\cdots\otimes X_{t_k}^\ast=\mathrm{sgn}(t)(Y^\ast\otimes X_{t_1}\otimes\cdots\otimes X_{t_k}),$$
where $Y\in R^{(0)}$, $X_{t_i}$ are the elements of the basis of the canonical $(-1)$-grading of $M_2$, $i=1,\ldots,k$, $t=t_1\cdots t_k\in T$, $\mathrm{sgn}(t)=\pm 1$, depending on the cases in Lemma \ref{L6}.
\epf 
\end{itemize}
\end{theorem}

\bigskip

As we have seen, not every grading of $L=\rsl(n)\subset R=M_n$ is induced by a grading of $R$. The best way to see why this happens in the ``A'' case but not in the case of ``B'',`` C'', ``D'' is to consider the connection between the gradings by an abelian group $G$ and the action by automorphism of the dual group $\widehat{G}$.

\section{Dual group action}\label{sdga}

 Let $F$ be a field. Denote by $\widehat G$ the dual group for $G$. Thus the elements
of  $\widehat G$ are all
irreducible characters $\chi:G\to F^*$, where $F^\ast$ is the
multiplicative group of the field $F$. 

If $\Lambda$ is a subset of $\wh{G}$ then we denote by $\Lambda^\perp$ the subgroup of $G$ given by
\bee{eperp}
\Lambda^\perp=\{ g\in G\,|\, \lambda(g)=1\mbox{ for all }\lambda\in\Lambda\}.
\ene
Similarly one defines $S^\perp$, a subgroup in $\wh{G}$, for any subset $S\subset G$. 

If $G$ is a a finite abelian group and $\wg$ is its group of characters over a field $F$ containing a primitive root of degree $n=|G|$ then $\wg\cong G$. More precisely,  if
$$
G=\langle a_1\rangle_{n_1}\times\cdots\times\langle a_j\rangle_{n_j}\times\cdots\times\langle a_k\rangle_{n_k}
$$ 
then
$$
\wh{G}=\langle \xi_1\rangle_{n_1}\times\cdots\times\langle \xi_i\rangle_{n_i}\times\cdots\times\langle \xi_k\rangle_{n_k}
$$ 
where $\xi_j(a_j)$ is a primitive root of $1$ of degree $n_j$ and $\xi_j(a_i)=1$ if $j\neq i$. 

Obviously, $\wg$ separates the elements of $G$ in the sense that for any $g_1\neq g_2$ there exists a $\chi\in\wg$ such that $\chi(g_1)\neq\chi(g_2)$. Then there is a natural isomorphism between $G$ and $\wh{\wg}$ given by $g\mapsto\psi_g$ where $\psi_g(\chi)=\chi(g)$.

Notice that if $\Lambda$ is a subgroup of $\wg$ and $H=\Lambda^\perp$ then there is a natural isomorphism $\alpha:\Lambda\rightarrow\wh{G/H}$ given by $\alpha(\lambda)(gH)=\lambda(g)$ where $\lambda\in\Lambda$. Its inverse $\beta:\wh{G/H}\rightarrow\Lambda$ is given by $\beta(\pi)(g)=\pi(gH)$. This, in particular, shows that $|\Lambda|\cdot|\Lambda^\perp|=|G|$.

We briefly recall the relation between the $G$-gradings and $\widehat
G$-actions in the case where $G$ is finite (see e.g. \cite{BSZ}) and $F$ contains a primitive root of degree $n=|G|$. Let $R$ be an arbitrary, i.e. not
necessarily associative, algebra graded by a finite abelian group $G$,
$R=\oplus_{g\in G}
R_g$. Any
element $a$ in $R$ can be uniquely decomposed as the
sum of homogeneous components, $a=\sum_{g\in G} a_g, a_g\in R_g$. Given
$\chi\in\widehat G$ we can define

\bee{e1}
\chi\ast a=\sum_{g\in G} \chi(g)a_g.
\ene

It is easy to observe that (\ref{e1}) defines a
$\widehat{G}$-action on $R$ by automorphisms and a subspace $V\subset
R$ is a
graded subspace if and only if $V$ is invariant under this action, i.e.
$\widehat{G}\ast V=V$. In particular, 
$a$ is a homogeneous iff $a$ is a common eigenvector for all $\chi\in\widehat{G}$.
The elements of the identity component $R_e$, $e$ the identity element
of $G$, are precisely the fixed points of the above action.

We will also use the following relation between actions by the subgroups of $\wh{G}$ and gradings
by the factor-groups of $G$. Notice that  thanks to the isomorphism  $\alpha:\Lambda\rightarrow\wh{G/H}$ we have the action of $\Lambda$ on $G/H$ given by $\lambda(gH)=\lambda(g)$ and on $R$ by $\lambda\ast(\sum_{x\in G/H}r_x)=\sum_{x\in G/H}\lambda(x)r_x$.

\begin{lemma}\label{F} Let $H$ and $\Lambda$ be subgroups of $G$ and $\wh{G}$, respectively, such that $H=\Lambda^\perp$. Then $a\in R$ is homogeneous in the natural $G/H$-grading if and only if $a$ is a common eigenvector for all $\chi\in\Lambda$. Similarly, $V\subset R$ is a $G/H$-graded subspace if and only if $\Lambda\ast V\subset V$. 
\end{lemma}

\pp  In view of what was said about the connection between gradings and actions, both claim will follow if we prove that the isomorphisms $\alpha:\Lambda\cong\wh{G/H}$ and $\beta:\wh{G/H}\cong\Lambda$, are compatible with the restriction of the action of $\wh{G}$ to  $\Lambda$ and the natural action of $\wh{G/H}$ on the $G/H$-graded algebra $R$. It is sufficient to check this compatibility for one of them, say, $\beta$. For instance, every element is the sum of those ones which are homogeneous in the $G/H$-grading. If $r$ has degree $gH$ then $r=\sum_{h\in H}r_{gh}$ and for $\pi\in \wh{G/H}$ we should have:
\begin{eqnarray*}
\pi\ast r&=&\pi(gH)r=\beta(\pi)(g)r=\sum_{h\in H}\beta(\pi)(g)r_{gh}\\&=&\sum_{h\in H}\beta(\pi)(gh) r_{gh}=\sum_{h\in H}\beta(\pi)\ast r_{gh}=\beta(\pi)\ast r.
\end{eqnarray*}
It follows that indeed the action of $\Lambda$ is equivalent to the action of the whole of $\wh{G/H}$, which is what we needed.
\epf

Now the difference between the ``A'' case and those of ``B'',`` C'', ``D'' in the case of an algebraically closed field $F$ of characteristic zero is that each algebra $L$ from the latter list, with the exception of some small rank algebras, has a canonical embedding in a matrix algebra $R$ in such a way that any automorphism of $L$ is induced by an automorphism of $R$. Our example above shows that this cannot be done in the case of simple Lie algebras of type ``A''. Still the following result from the classical Lie Theory tells us the following.

\begin{theorem}{\rm (\cite[Chapter IX, Theorem 5]{NJL})} 
The automorphism group of any Lie algebra $L=\rsl(n)\subset M_n$ over an algebraically closed field $F$ of characteristic zero is generated by the inner automorphisms $X\rightarrow T^{-1}XT$, $T$ a nondegenerate matrix in $M_n$ and an outer automorphism of order 2, $X\rightarrow - ^t\! X$, where $\tr{X}$ is the transpose of $X$. In the case $n=2$ this latter mapping of $L$ is also induced by an inner automorphism of $M_2$.
\end{theorem}

Now suppose we are given a grading of $L=\rsl(n)$ by a finite abelian group $G$. We consider the dual group $\wh{G}$. Then we have the action  of $\wh{G}$ by the automorphisms of $L$. We will call a $G$-grading $L=\bigoplus_{g\in G}L_g$ on $L=\rsl(n)$ \textit{inner} if $\wg$ acts on $L$ by inner automorphisms. Otherwise we call this grading \textit{outer}.

Suppose first that the grading is inner, that is, for each $\chi\in  \wh{G}$ there is a nondegenerate matrix $T_\chi\in M_n$ such that

\bee{nf1}
\chi\ast X=T_\chi^{-1}XT_\chi,\mbox{  for any }X\in L.
\ene
 
It is then obvious that the same formula (\ref{nf1}) defines an action of $\wh{G}$ on $R$, thus a $G$-grading of $R=M_n$. This will be a unique $G$-grading of $R=M_n$ that induces the given grading of $L=\rsl(n)$, and it is given by $R_e=FI\oplus L_e$ and $R_g=L_g$ for $g\neq e$. The description of such gradings of $\rsl(n)$ is identical to that of $M_n$ given in Theorem \ref{tbsz}. For convenience, we formulate this result below in the language of Lie algebras.

\begin{theorem}\label{tbszl}Let $F$ be an algebraically closed field of characteristic zero, $G$ a finite abelian group. If $L=\rsl(n)$ and $L=\bigoplus_{g\in G}L_g$ is an inner $G$-grading then there is a grading  $R=\bigoplus_{g\in G}R_g$ of $R=M_n$ described in Theorem \ref{tbsz} such that $L_g=R_g\cap L$ for any $g\in G$.\epf
\end{theorem}

\section{Outer gradings. General Results}

Main results and arguments of the remaining sections of this paper heavily depend on \cite{BShZ}, where it is assumed that the base field $F$ is algebraically closed of characteristic zero. Therefore, although some intermediate results and arguments could be proved under milder restrictions, we still prefer from now on to work under this assumption. 

So we have to consider the case of an outer grading, namely where there is an element of $\wh{G}$ which acts on $L$ as an outer automorphism. In this case there is a subgroup $\Lambda\subset \wh{G}$ of index 2, which acts by inner automorphisms on $L$. An element, say $\vp\in \wh{G}\setminus \Lambda$ is the composition of an inner automorphism $X\rightarrow \Phi^{-1} X\Phi$ and the canonical automorphism $X\rightarrow - ^t\! X$. Thus the action of $\vp$ is given by

\bee{nf2} 
\vp\ast X= -\Phi^{-1}( \tr{X})\Phi.
\ene

This formula defines also an action of $\vp$ on the associative algebra $R=M_n$ by a Lie automorphism such that $\vp\ast I = - I$.
Moreover, we now can extend the action of $\wh{G}$ on $R$ by Lie automorphisms if we set $\chi\ast I=I$ for $\chi\in \Lambda$ and $\chi\ast I=-I$ otherwise. 

Notice that in (\ref{nf2}) the matrix $\Phi$ is defined up to a scalar multiple. A simple but useful remark about the form of $\Phi$ is the following.
\begin{lemma}\label{nl1}
If one applies an inner automorphism to $M_n$ induced by a matrix $C$ then $\Phi$ in (\ref{nf2}) is changed to
$\tr{C}\Phi C$. In other words, $\Phi$ behaves as the matrix of a bilinear form.
\end{lemma}
\pp Indeed, we must have $\vp\ast\widetilde{X}=\widetilde{\vp\ast X}$ where tildes mean the matrices modified as a result of the inner automorphism in question. If we denote by $\Psi$ the modified $\Phi$ then we will have
$$-\cj{\Psi}{\tr{(\cj{C}{X})}}=-\cj{C}{\cj{\Phi}{\tr{X}}}\mbox{ or } \cj{(\tr{C^{-1}}\Psi)}{\tr{X}}=\cj{(\Phi C)}{\tr{X}}.$$
It follows that $\tr{(C^{-1})}\Psi=\alpha \Phi C$, for some coefficient $\alpha$ and so we may take $\Psi =\, \tr{C}\Phi C$, as claimed.\bx

\begin{proposition}\label{np1} One can choose the character $\vp$, the elements $a,h\in G$, the subgroups $\Lambda_1\subset\wh{G}$ and $K\subset G$ in such a way that the following hold:
\begin{enumerate}
\item[{\rm (i)}] The order of $\vp$ is a 2-power: $o(\vp)=2m=2^k$ for a natural $k\geq 1$;
\item[{\rm (ii)}] $\wh{G}=\langle\vp\rangle\times\Lambda_1$ and $\Lambda=\langle\vp^2\rangle\times \Lambda_1$;
\item[{\rm (iii)}] $G=\langle a\rangle\times K$ so that $\vp(K)=1$, $\vp(a)=\rho$, $\rho$ a $(2m)^\mathrm{th}$ primitive root of $1$, $\Lambda^\perp=\langle h\rangle$ where $h=a^{m}$, $\Lambda_1^\perp=\langle a\rangle$, $\vp(h)=-1$.
\end{enumerate}
\end{proposition}
\pp This is an easy exercise, which we accomplish here for completeness. Let us denote $\wh{G}$ by $\Gamma$. We can decompose $\Gamma$ as $\Gamma=\Gamma_2\times\Gamma_2^\prime$ where $\Gamma_2$ is the Sylow $2$-subgroup and $\Gamma_2^\prime$ is the product of all other Sylow $p$-subgroups. Since $\Gamma^2\subset\Lambda$ it easily follows that $\Gamma_2^\prime\subset\Lambda$. Now we can decompose $\Gamma_2$ as the direct product of cyclic subgroups
\bee{fag}
\Gamma_2=\langle \chi_1\rangle_{k_1}\times\cdots\times\langle \chi_i\rangle_{k_i}\times\cdots\times\langle \chi_m\rangle_{k_m}\mbox{ with }k_1\,|\ldots|\,k_i\,|\ldots|\, k_m
\ene
where $i$ is the smallest index such that $\chi_i\notin \Lambda$. If we replace in (\ref{fag}) each $\chi_j\notin \Lambda$ by $\chi_i^{-1}\chi_j$ then the decomposition remains but only one direct summand will be outside $\Lambda$. Now we can set $\vp=\chi_i$,  and choose $\Lambda_1$ to be the direct product of the remaining cyclic factors of the transformed $\Gamma_2$, and $\Gamma_2^\prime$. This then proves (i) and (ii).

For the proof of (iii) the easiest way is to recall that $\wh{\wh{G}}\cong G$. Then if we have a direct cyclic decomposition 
$$
\wh{G}=\langle \xi_1\rangle_{n_1}\times\cdots\times\langle \xi_i\rangle_{n_i}\times\cdots\times\langle \xi_k\rangle_{n_k}
$$ 
 of $\wh{G}$ it causes a similar decomposition 

$$
G\cong\wh{\wh{G}}=\langle a_1\rangle_{n_1}\times\cdots\times\langle a_j\rangle_{n_j}\times\cdots\times\langle a_k\rangle_{n_k}
$$ 
of $G\cong\wh{\wh{G}}$ such that $\xi_j(a_j)$ is a primitive root of $1$ of degree $n_j$ and $\xi_j(a_i)=1$ if $j\neq i$. Since we can assume $\varphi$ being
    one of $\xi_i$, our claim (iii) is an immediate consequence of these two decompositions.\bx

Let us denote by $H$ the subgroup of order 2 generated by $h$, as in the preceding Proposition \ref{np1}. Then $\Lambda\cong \wh{G/H}$, with a well-defined isomorphism $\alpha:\Lambda\rightarrow\wh{G/H}$ given by $\alpha(\lambda)(\bar{g})=\lambda(g)$. Here $\bar{g}=gH$, an element of $G/H$. Therefore it is easy to check that we have a $G/H$-grading on $L=\rsl(n)$ defined by $L_{\bar{g}}=L_g\oplus L_{gh}$. As noted above, we have an extension of the action of $\wh{G}$ on $R=M_n$ by Lie automorphisms if we set $\chi\ast I=I$ for $\chi\in\Lambda$ and $\chi\ast I=-I$ otherwise. This is equivalent to saying that the outer grading of $L$ extends to a Lie grading of $R$ if we set $R_h=FI\oplus L_h$ and $R_g=L_g$ if $g\neq h$. For the factor-grading we then have $R_{\bar{e}}=FI\oplus L_{\bar{e}}$ and $R_{\bar{g}}=L_{\bar{g}}$ if $\bar{g}\neq \bar{e}$.

\begin{lemma}\label{lelfin}
Let $\rsl(n)=L=\bigoplus_{g\in G}L_g$ be an outer grading on $L$ and $G,H,\Lambda,\vp$ be as in Proposition \ref{np1}. Then the $G/H$-grading on $L$ is inner and induced from a $G/H$-grading on $R=M_n$. Moreover, if $R=R^{(0)}\otimes R^{(1)}\otimes\cd\otimes R^{(k)}$ as in Theorem \ref{tbsz}, with respect to $G/H$-grading, then all $L\cap R^{(i)}$, $i=0,1,\ld,k$, and $L\cap (R^{(1)}\otimes\cd\otimes R^{(k)})$ are $G$-graded subalgebras of $L$.
\end{lemma}

\pp Since $\Lambda\cong \wh{G/H}$ acts on $L$ by inner automorphisms it follows that the $G/H$-grading of $L$ is inner hence induced from a unique $G/H$-grading of $R$. It follows that the $G/H$-Lie-grading of $R$ described just before this lemma is actually an associative grading. If we decompose $R$ as the tensor product $R=R^{(0)}\otimes R^{(1)}\otimes\cd\otimes R^{(k)}$ following Theorem \ref{tbsz} then all $R^{(i)}$ are stable under $\Lambda$-action. Since $\Lambda$ and $\vp$ commute, and each $R_{\bar{g}}$ is a weight subspace under the action of $\Lambda$, it follows that for each $\bg\in G/H$ we have $\vp\ast R_{\bg}\subset R_{\bg}$.  Moreover, 
$$
R^{(0)}=\bigoplus_{\bg\,\in\,\su {R^{(0)}}}R_{\bg}\mbox{ hence }\vp\ast R^{(0)}=R^{(0)}.
$$
Since $\su {R^{(i)}}\cap \su {R^{(j)}}=\{ \be\}$ for any $i\neq j$  and the centralizer $C_R(R^{(0)})$ of $R^{(0)}$ in $R$ equals  $R^{(1)}\otimes\cd\otimes R^{(k)}$ it follows that if $R^{(i)}_{\bg}\neq 0$, $i\geq 1$, and $\bg\neq \be$ then $R^{(i)}_{\bg}=R_{\bg}\cap C_R(R^{(0)})$. It is immediate then that $\vp\ast R^{(i)}_{\bg}=R^{(i)}_{\bg}$ for all $\bg\neq\be$ and $1\leq i\leq k$. Now let us consider $R^{(i)}_{\be}$. If $i=1,\ld,k$ then $R^{(i)}_{\be}=\Sp{I}$. Since $\vp\ast I=-I$ we have  $\vp\ast R^{(i)}_{\be}=R^{(i)}_{\be}$ for any $i=1,\ld,k$. Finally, from the intersection property of the supports it follows that 
$$
R_{\be}=R^{(0)}_{\be}\otimes R^{(1)}_{\be}\otimes\cd\otimes R^{(k)}_{\be}=R^{(0)}_{\be}.
$$
As before, we then have $\vp\ast R^{(0)}_{\be}=R^{(0)}_{\be}$. Now if $S=R^{(1)}\otimes\cd\otimes R^{(k)}$ then it follows from $\vp(xy)=-\vp(y)\vp(x)$ that $\vp\ast S=S$. As a result, all $R^{(i)}$ and $R^{(1)}\otimes\cd\otimes R^{(k)}$ are stable under $\wh{G}$.

Finally, since the action of $\wg$ leaves $L$ invariant we have that each $L\cap  R^{(i)}$, $i=0,1,\ld,k$, and $L\cap (R^{(1)}\otimes\cd\otimes R^{(k)})$ are $G$-graded subalgebras of $L$, as required.\epf

Since the action of $\Lambda$ on $L$ is inner with each $\lambda\in \Lambda$ one can associate a non-degenerate matrix $T_\lambda\in M_n$ such that \bee{nf3}\lambda\ast X=\cj{T_\lambda}{X},\mbox{ for any }X\in L.\ene We can also say that $X$ is homogeneous of degree $\bar{g}$ if and only if

 \bee{nf4}
\cj{T_\lambda}{X}=\lambda(g)X.
\ene

Now, as before, we can extend the action of $\Lambda$ on $L$ to $R=M_n$. It is given by the same formula (\ref{nf3}). Therefore, we have a $G/H$-grading of $M_n$. Actually, as mentioned before, for any $\bar{g}\neq\bar{e}$ we have that $R_{\bar{g}}=L_{\bar{g}}$ while $R_{\bar{e}}= FI\oplus L_{\bar{e}}$. 

If we know the $\bar{G}=G/H$-grading of $L$ and an automorphism $\vp$ satisfying $\vp\ast L_{\bar{g}}=L_{\bar{g}}$ then we can recover the $G$-grading using the following procedure. We recall that $L_{\bar{g}}=L_g\oplus L_{gh}$. 

\begin{proposition}\label{p00} Using our previous notation, one can find $L_a$, for any $a\in \bar{g}$, in the form as follows:
\bee{nf8}
L_a=\{ X+\vp(a)^{-1}(\vp\ast X)\,\vert\, X\in L_{\bar{g}}\}.
\ene
\end{proposition}

\pp Indeed, let $M_a$ denote the right side of (\ref{nf8}). Since any element in $M_a$ is still in $L_{\bar{g}}$ the action of any $\lambda\in\Lambda$  amounts to the scalar multiplication by $\lambda(g)=\lambda(a)$, for any $a\in \bar{g}$. Now let us check the action of $\vp$.

$$
\vp\ast (X+\vp(a)^{-1}(\vp\ast X))=\vp\ast X+\vp(a)^{-1}(\vp^2\ast X)).
$$

Since $\vp^2\in \Lambda$ we have that $\vp^2\ast X=\vp^2(\bg)X=\vp^2(a)X=(\vp(a))^2X$. Plugging this value in (\ref{nf8}) we obtain
$$
\vp\ast (X+\vp(a)^{-1}(\vp\ast X))=\vp(a)(\vp(a)^{-1}(\vp\ast X)+X).
$$
Since $\vp$ and $\Lambda$ generate the whole of $\wh{G}$ it follows that $X+\vp(a)^{-1}(\vp\ast X)$ is an eigenvector for the action of any $\chi\in\wh{G}$ with eigenvalue $\chi(a)$, proving that the degree of this element in the $G$-grading is $a$. So $M_a\subset L_a$. Notice also that
$$
M_{ah}=\{ X+\vp(ah)^{-1}(\vp\ast X)\,\vert\, X\in L_{\bar{g}}\}=\{ X-\vp(a)^{-1}(\vp\ast X)\,\vert\, X\in L_{\bar{g}}\}.
$$

It is immediate then that $M_a+ M_{ah}=L_{\bar{g}}$. Since, being graded components of the $G$-grading of $L$,  the subspaces $L_a$ and $L_{ah}$ have trivial intersection we must conclude that $M_a=L_a$ and $M_{ah}=L_{ah}$, so that our claim is correct.\epf

Now we know that the action of $\vp$ on $L$ is given by (\ref{nf2}) and that both $\vp$ and $\Lambda$ belong to the same abelian group $\wh{G}$. This causes a number of relations between the matrix $\Phi$, as in (\ref{nf2}) and the matrices $T_\lambda$, as in (\ref{nf3}), for each $\lambda\in\Lambda$.

We have
$$\lambda\ast(\vp\ast X)=\lambda\ast(-\cj{\Phi}{\tr{X}})=-\cj{T_\lambda}{\cj{\Phi}{\tr{X}}}=\cj{(\Phi T_\lambda)}{\tr{X}}.$$
and also
$$\vp\ast(\lambda\ast X)=\vp\ast(\cj{T_\lambda}{X})=-\cj{\Phi}{\tr{(\cj{T_\lambda}{X})}}=\cj{(\tr{\iv{T}}\Phi)}{\tr{X}}.$$
Since $X$ is any matrix with trace zero, it follows that the conjugating matrices on the right sides of both equations differ only by a scalar $\beta$. We then have

\bee{nf5}
 \Phi T_\lambda=\beta\tr{\iv{T_\lambda}}\Phi\;\mbox{ or }\;\tr{T_\lambda}\Phi T_\lambda=\beta\Phi.
\ene

Now let us extend the action of $\vp$ to $M_n$, using the same formula (\ref{nf2}). It is obvious that $\vp$ and $\Lambda$ remain commutative. It follows then that for any homogeneous component $L_{\bar{g}}$, respectively, $R_{\bar{g}}$  we have  $\vp\ast L_{\bar{g}}=L_{\bar{g}}$, respectively,  $\vp\ast R_{\bar{g}}=R_{\bar{g}}$. Let us notice that the action of $\vp$ on $M_n$ satisfies the following relation: 

\bee{nf6}
 \vp(XY)=-\vp(Y)\vp(X)\mbox{  for all  }X,Y\in M_n.
\ene

It is clear, conversely, that any mapping $\vp:M_n\rightarrow M_n$ satisfying (\ref{nf6}) must be of the form (\ref{nf2}). Indeed, if we combine such $\vp$ with $X\rightarrow -\tr{X}$ then we obtain an automorphism of $M_n$, hence a conjugation by an appropriate non-degenerate matrix $\Phi$, as needed.

If $\vp$ is of order $2$ then $-\vp$ becomes an involution and  so the factor-grading on $M_n$ by $G/H$ becomes an involution grading. Such gradings have been completely described in Subsection \ref{sAIGMA}.

We cannot immediately apply these results to outer gradings of $L=\rsl(n)$ since our outer automorphism $\vp$ need not be an element of order $2$. But there is a way to replace, in certain situations, outer automorphisms of arbitrary order by those of order two. Notice that in the statement of the theorem that follows we model on the form of $G$ and $\wg$ as it is given in Proposition \ref{np1}.

\begin{theorem}\label{t00} Let $L=\rsl(n)\subset R=M_n$ be graded by a group $G$ of the form $G=\langle a\rangle_{2m}\times K$. Suppose the dual group has the form of $\widehat{G}=\langle\varphi\rangle_{2m}\times\Lambda_1$ where $\varphi(a)=\rho$, a primitive root of degree $2m$, $\varphi(K)=1$, $\Lambda_1(a)=1$, $\Lambda_1\cong\widehat{K}$. Let the action of $\varphi$ on $L$ be outer, while that of $\Lambda=\langle \vp^2\rangle\times \Lambda_1$ inner. Suppose $\psi$ is an inner automorphism of $L$ such that the action of $\vp^2$ coincides with $\psi^2$ and $\psi$ commutes with the action of $\wg$. Then if we set $h=a^m$, there is an involution $\ast$ of $R$ and an involution preserving $G$-grading 
$$
R=\sum_{g\in G}\widetilde{R}_g
$$
of $R$ such that 
\bee{fmain}
L_g=\widetilde{R}_g^{(-)}\oplus(\widetilde{R}_{gh}^{(+)}\cap L).
\ene
for all $g\in G$.
\end{theorem}

\pp Note that by Corollary \ref{cout} the above formulas define a $G$-grading of $L$.

To prove the converse, we consider a new abelian group $\widetilde{G}=\langle c\rangle_2 \times G$ and its dual $\widehat{\widetilde{G}}=\langle \omega\rangle_2 \times \wg$. The action of $\widehat{\widetilde{G}}$ on $\widetilde{G}$ naturally extends that of $\widehat{G}$ on $G$ by setting $\omega(c)=-1$.

Now the action of $\widehat{G}$ on $L$ naturally extends to that of $\widehat{\widetilde{G}}$ if we set $\omega\ast X=(\varphi\ast\psi^{-1})(X)$. As a consequence, $L$ acquires a $\widetilde{G}$ grading given by
\begin{equation}\label{ecg}
L_{(c^i,g)}=\{ X\in L\,\vert\,(\omega^j,\chi)\ast X=(-1)^{ij}\chi(g)X\}.
\end{equation}
Here 
$$
(\omega,\chi)\ast X=\omega\ast(\chi\ast X)=\chi\ast(\omega\ast X)\mbox{ and }i,j=0,1,\;\chi\in \widehat{G}.
$$
It is obvious that then if we know the $\widetilde{G}$-grading of $L$ then the $G$-grading can be recovered by setting
$$
L_g=L_{(e,g)}\oplus L_{(c,g)}.
$$

To find the components of the $\widetilde{G}$-grading of $L$, we proceed as follows. We rewrite $\widehat{\widetilde{G}}$ in the form 
$$
\widehat{\widetilde{G}}=\langle \omega\rangle_2 \times \langle \vp\omega^{-1}\rangle_{2m}\times \Lambda_1.
$$
In this case, the subgroup $P=\langle \vp\omega^{-1}\rangle_{2m}\times \Lambda_1$ acts on $L$ by inner automorphisms. As previously, this action naturally extends to an associative action on $R$. Under this extended action $-\omega$ is an involution $\ast$ of an associative algebra $R$ and so we write $X^\ast=-\omega\ast X$. The components of the $\widetilde{G}$-grading of $L$ are the symmetric and skew-symmetric components of the components of the factor-grading by  $\widetilde{G}/P^\perp$, under this involution (see Lemma \ref{F}).

Now
$$
P^\perp=\langle \varphi\omega^{-1}\rangle^\perp\cap(\Lambda_1)^\perp.
$$
Since $(\Lambda_1)^\perp=\langle c\rangle_2\times\langle a\rangle_{2m}$ and $\omega^{-1}=\omega$, we need only the annihilator of $\langle \varphi\omega\rangle$ on $\langle c\rangle_2\times\langle a\rangle_{2m}$, which is, obviously, $$Q=\langle (c,a^m)\rangle_2=\langle (c,h)\rangle_2.$$ So we have to determine the $\widetilde{G}/Q$-grading of $L$ and then partition it into the skew-symmetric and symmetric components. Obviously, $\widetilde{G}/Q\cong G$, so that actually we have an inner $G$-grading. 

We consider the components of the $\widetilde{G}/Q$-grading of $L$ by $L_{(c^i,g)Q}$. Then
$$
L_{(c^i,g)Q}=L_{(c^i,g)}\oplus L_{(c^{i+1},\,gh)}\mbox{ where }i=0,1\,(\mbox{mod}\,2),\; g\in G.
$$
Now we can apply Proposition \ref{p00} to conclude that for any $X\in L_{(e,g)Q}$ the element $$X+\omega((e,g))^{-1}(\omega\ast X)$$ will be of degree $(e,g)$ while $$X+\omega((c,gh))^{-1}(\omega\ast X)$$ of degree $(c,gh)$. Therefore, the elements of the form $X+\omega\ast X$, that is the skew-symmetric elements of $\ast=-\omega$ are of degree $(e,g)$ while the elements of the form $X-\omega\ast X$, that is the symmetric elements are of degree $(c,gh)$.

Now $L_g=L_{(e,g)}\oplus L_{(c,g)}$. By the above, $L_{(e,g)}$ is the set of skew-symmetric elements in $L_{(e,g)Q}$ while $L_{(c,g)}$ consists of the symmetric elements in $L_{(c,g)Q}$, which is the same as $L_{(e,gh)Q}$.

As we mentioned above, $\widetilde{G}/Q\cong G$ allows us to make the above inner grading of $L$ by $\widetilde{G}/Q$ into an inner grading by $G$. We may set $\widetilde{L}_g=L_{(e,g)Q}=L_{(c,gh)Q}$. Because this grading of $L$ is inner it is the restriction of the $G$-grading of $R$ defined by $\widetilde{R}_e=FI\oplus \widetilde{L}_e$ and $\widetilde{R}_g=\widetilde{L}_g$ for $g\neq e$. Our involution $-\omega$ of $R$ is then compatible with this new grading. If we denote by $\widetilde{R}^{(\pm)}_g$ the set of $-\omega$-symmetric (resp., skew-symmetric) elements of $\widetilde{R}_g$, and recall that skew-symmetric elements of an involution always have trace 0, we will have
$$
L_g=\widetilde{R}_g^{(-)}\oplus(\widetilde{R}_{gh}^{(+)}\cap L).
$$
\epf

\begin{corollary}\label{c00} Suppose the outer automorphism $\vp$ as above is of order $2$. Then $-\vp$ is an involution $\ast$ and there is an inner involution compatible $G$-grading  
$$
R=\sum_{g\in G}\widetilde{R}_g
$$
of $R$ such that 
\bee{e0005}
L_g=\widetilde{R}_g^{(-)}\mbox{ if }g\in K\mbox{ and }L_g=\widetilde{R}_{gh}^{(+)}\cap L\mbox{ if }g\notin K.
\ene
\end{corollary}

\pp In this case we have $\psi=\mathrm{id}$, $\omega$ acts as $\vp$, and $P=\langle\vp\omega\rangle\times\Lambda_1$. Now $\widetilde{L}_g=L_{(e,g)Q}=L_{(e,g)}\oplus L_{(c,gh)}$. By definition (\ref{ecg}), $L_{(e,g)}=\{ X\in L_g\vert \omega^j\ast X=X\}$, for any $j=0,1$. Now $\omega$ acts as $\vp$ and so $L_{(e,g)}=\{ X\in L_g\vert -X=X\}=0$ for $g\notin K$. Also, $L_{(e,g)}=L_g$ for $g\in K$. Similarly, $L_{(c,g)}=\{ X\in L_g\vert \omega^j\ast X=(-1)^jX\}$. Now 
$L_{(c,g)}=\{ X\in L_g\vert X=-X\}=0$ for $g\in K$. Also, $L_{(c,g)}=L_g$ for $g\notin K$. By the proof of our theorem, $L_{(e,g)}$ is the set of skew-symmetric elements in $\tilde{L}_g$ while $L_{(c,g)}$ is the set of symmetric elements in $\widetilde{L}_{gh}=\widetilde{R}_{gh}\cap\widetilde{L}$ and our lemma follows.
\epf

To study outer gradings on $L=\rsl(n)$ we are going to apply Lemma \ref{lelfin}. According to this lemma, the matrix algebra $R=M_n$ possesses a grading by the group $\bar{G}=G/H$. One can decompose $R$ as the tensor product $R=P\otimes Q$ of $\bar{G}$-graded subalgebras $P=R^{(0)}\cong M_p$ and $Q=R^{(1)}\otimes\cdots\otimes R^{(k)}\cong M_q$, the $\bar{G}$-grading on $P$ being elementary and the $\bar{G}$-grading on $Q$ being fine. Moreover, both $P$ and $Q$ are invariant under $\vp$. The intersections $L\cap P$ and $L\cap Q$ are thus $G$-graded Lie algebras $\rsl(p)\subset P=M_p$ and $\rsl(q)\subset Q=M_q$, respectively, and the $G/H$-grading on $P$ is inner and elementary while $\bar{G}$-grading on $Q$ is inner and fine. Since the action of $\vp$ satisfies $\vp(XY)=-\vp(Y)\vp(X)$ on $R$, it follows that on each $P$ and $Q$ we can write $\vp$ in the form (\ref{nf2}) for appropriate $\Phi$.

These remarks allow us to restrict ourselves to the cases where the original $G$-grading is either such that the induced $G/H$-grading is elementary or such that the induced $G/H$-grading is fine. When we get information about these two cases we will have to return to the general case and consider the tensor products. 

\section{Fine gradings}\label{s2}

We start with a lemma close to \cite[Lemma 4]{BShZ}.

\begin{lemma}\label{fin1}
Let $R=M_n(F)=\bigoplus_{t\in T}R_t$ be the $n\times n$-matrix algebra with an $\ve$-grading,
$T=\langle a\rangle_n\times \langle b\rangle_n$. Let also $\vp: R\rightarrow R$ be a mapping on $R$ defined
by $\vp*X=-\Phi^{-1~t}X\Phi$. If $\vp*R_t=R_t$ for all $t\in T$ then $n=2$ and $\vp$ acts
on $\rsl(2)$ as the conjugation  by one of the matrices $X_a$, $X_b$ or $X_{ab}$ (see (\ref{e00})).
\end{lemma}

\pp First we consider the $\vp$-action on $X_a$. Since $R_a$ is stable under $\vp$,
$$
-\Phi^{-1~t}X_a\Phi=-\Phi^{-1}X_a\Phi=\alpha X_a
$$
for some scalar $\alpha\ne 0$. Then
\bee{f1}
X_a^{-1}\Phi X_a=\beta\Phi
\ene
with some
$\beta=-\alpha^{-1}$. Since $X_a^n=I$, we obtain $\beta^n=1$, so that $\beta=\ve^j$ for some $0\le j\le n-1$.

Denote by $P$ the linear span of $I,X_a,\ld,X_a^{n-1}$. Then $R=P\oplus X_bP\oplus \cd
\oplus X_b^{n-1}P$ as a vector space and the conjugation by $X_a$ acts on $X_b^iP$ as the
multiplication by $\ve^{-i}$. In particular, all eigenvectors with eigenvalue $\ve^{-j}$
are in $X_b^jP$. It follows that $\Phi\in X_b^jP$, that is, $\Phi= X_b^jQ$ for some $Q\in P$.

Now we consider the action of $\vp$ on $X_b$:
$$
\vp* X_b=-\Phi^{-1~t}X_b\Phi=-\Phi^{-1}X_b^{-1}\Phi=\gamma X_b,
$$
that is, $X_b\Phi X_b = \mu\Phi$ with $\mu=-\gamma^{-1}\ne 0$. If we write $Q=\sum\alpha_i
X_a^i$ then
\bee{f2}
X_b\Phi X_b = X_b^j\sum_{i}\alpha_iX_b X_a^i X_b =
X_b^j\sum_{i}\alpha_i' X_a^i X_b^2=\mu\Phi=\mu X_b^j\sum_{i}\alpha_iX_a^i,
\ene
In this case $X_b^j\sum_{i}\alpha_i' X_a^i X_b^2=\mu X_b^j\sum_{j}\alpha_iX_a^i$
where the scalars $\alpha_j'$ can be explicitly
computed using (\ref{a1}). Since the degrees in $X_a,X_b$ define the degrees in the $T$-grading, we can
see that (\ref{f2}) immediately implies $X_b^2=I$, i.e. $n=2$.

As we have shown before, (\ref{f1}) implies $\Phi=X_b^j Q$ with $Q=\alpha_0
I+\alpha_1 X_a$. Since $n=2$, the argument following (\ref{f1}) applies if we change places of $a$ and $b$  so that $\Phi=X_a^k(\beta_0 I +
\beta_1 X_b)$. Comparing these two expressions we obtain that 
$\Phi$ must be one of $I$, $X_a$, $X_b$, or $X_{ab}$, up to a scalar multiple. Finally note that
$X_{ab}^{-1}YX_{ab}=-^tY$ for any traceless $2\times 2$ matrix $Y$ and then
$-\Phi^{-1~t}Y\Phi$ is the conjugation by one of the matrices $X_{ab}\cdot I$,
$X_{ab}\cdot X_a$, $X_{ab}\cdot X_b$ or $X_{ab}^2$ which are scalar multiples of $X_{ab},
X_a,X_b$ and $I$,  respectively. \epf

A quick corollary is as follows.

\begin{corollary}\label{fin2}
Let $R=R^{(0)}\otimes R^{(1)}\otimes\cd\otimes R^{(k)}$ be the decomposition of
$R=M_n(F)$ corresponding to the factor-grading of $R$ by $\bG=G/H$ where
$R^{(0)}=M_{n_0}(F)$ has an elementary $\bG$-grading while the $\bG$-grading on $R^{(1)}\otimes\cd\otimes R^{(k)}$ is fine, as in Lemma \ref{lelfin}. Then $R^{(1)}\cong\ld\cong
R^{(k)}\cong M_2(F)$ have fine $(-1)$-grading. Also, the restriction of the outer
automorphism $\vp\in\hat G \setminus\Lambda$ to any factor $R^{(i)}, 1\le i\le k$,
coincides with the action of some inner automorphism. The order of the restriction of $\vp$ to $R^{(1)}\otimes\cd\otimes R^{(k)}$ equals $2$.
\end{corollary}

\pp Only the claim about the order requires some justification. To do this, we have to apply Lemma \ref{fin1} and the formula (\ref{nf6}). \epf

To formulate the next theorem we recall that if a grading by a group $G$ on a matrix algebra $R=M_n$ is fine then it follows from the results of Subsection \ref{sAGMA} that the support of the grading is a subgroup, say $T$, and a canonical graded basis can be chosen of nondegenerate matrices $X_t$, $t\in T$, such that $X_tX_{t'}=\alpha(t,t')X_{tt'}$ for a $2$-cocycle $\alpha:T\times T\rightarrow K^\ast$ . If there is an involution $\ast$ on $R$, which respects the grading, then there is a function $\beta: T\rightarrow \{\pm 1\}$ such that $X_t^\ast=\beta(t)X_t$. 

\begin{theorem}\label{tfin} Let $L=\rsl(n)=\bigoplus_{g\in G}L_g\subset R= M_n$ be given an outer grading by a finite abelian group $G$ such that the respective $\bG$-grading of $R$ is fine. Then $n=2^k$ and there is a fine involution grading on $R=\bigoplus_{g\in G}\widetilde{R}_g$ with a subgroup $T$ as its support and an element $h$ of order $2$ in $G$ such that $R$,  as a $G$-graded algebra with involution is isomorphic to the tensor product $R=R^{(1)}\otimes\cdots\otimes R^{(k)}$ of graded involution stable subalgebras $R^{(i)}$ each of which is a matrix algebra of order $2$ of one of the four types in Lemma \ref{L6}, with support $T_i\cong \mathbb{Z}_2\times \mathbb{Z}_2$. The support of $R$ is $T=T_1\times\cdots\times T_k$. The basis of $R$ is formed by the Kronecker products $X_t=X_{t_1}\otimes\cdots\otimes X_{t_k}$ where $t_i\in T_i$. Further, 
$$X_t^\ast = X_{t_1}^\ast\otimes\cdots\otimes X_{t_k}^\ast=\beta(t)X_{t_1}\otimes\cdots\otimes X_{t_k}\mbox{ and } \beta(t)=\beta(t_1)\cdots\beta(t_k).$$

Finally, the original $G$-grading of $L$ can be recovered as follows (we mention only nonzero components).

\emph{Case 1}: $h\notin T$. 
\begin{enumerate}
\item[$\mathrm{(a)}$] $L_t=\Sp{X_t}$ for $t\in T$ such that $\beta(t)=-1$;
\item[$\mathrm{(b)}$] $L_{g}=\Sp{X_t}$ for $g\in G\setminus\{ h\}$ such that $g=th$, $t\in T$ such that and $\beta(t)=1$.
\end{enumerate}

\emph{Case 2}: $h\in T$. 
\begin{enumerate}
\item[$\mathrm{(a)}$] $L_t=\Sp{X_t,X_{th}}$ for $t\in T\setminus\{ h\}$ such that $\beta(t)=-1$ and $\beta(th)=1$;
\item[$\mathrm{(b)}$] $L_t=\Sp{X_t}$ for $t\in T$ such that $\beta(t)=-1$ and $\beta(th)=-1$;
\item[$\mathrm{(c)}$] $L_t=\Sp{X_{th}}$ for $t\in T\setminus\{ h\}$ such that $\beta(t)=1$ and $\beta(th)=1$
\end{enumerate}
\end{theorem}

\pp By the argument preceding Lemma \ref{lelfin} and the definition of the factor-grading, $$R_{\be}=\Sp{I}\oplus L_{\be}=\Sp{I}\oplus L_e\oplus L_h.$$
By our hypotheses, $\dim R_{\be} =1$ and so $L_e=L_h=\{ 0\}$. Now we adopt the notation and the argument of Theorem \ref{t00}. By Corollary \ref{fin2}  any character $\vp\in G\setminus\Lambda$ acts as an automorphism of order $2$ and so Theorem \ref{t00} and Corollary \ref{c00} state that then $R$ possesses an inner grading $R=\oplus_{g\in G}\widetilde{R}_g$ stable under the involution $\ast = -\vp$. Following the proof of Theorem \ref{t00}, we write $G=\langle a\rangle_{2m}\otimes K$, its dual $\wh{G}=\langle\vp\rangle_{2m}\times\Lambda_1$, $\widetilde{G}=\langle c\rangle_2 \times G$ and its dual $\widehat{\widetilde{G}}=\langle \omega\rangle_2 \times \wg$. One chooses
$P=\langle \vp\omega^{-1}\rangle_{2m}\times \Lambda_1$ and $Q=P^\perp=\langle (c,a^m)\rangle_2=\langle (c,h)\rangle_2.$ Then $\widetilde{R}_g$ is defined by $\widetilde{R}_g=R_{(e,g)Q}$. Since $\vp$ commutes with $\widehat{\widetilde{G}}$ this grading is stable under the involution. To check that this grading of $R$ is fine we only need to check $\dim\widetilde{R}_e=\dim R_{(e,e)Q}=1$ (\cite[Section 4]{BSZ}). Now 
$$\Sp{I}\subset R_{(e,e)Q}=R_{(e,e)}\oplus R_{(c,h)}\subset R_e\oplus R_h=L_e\oplus L_h\oplus\Sp{I}=\Sp{I},$$
proving that indeed $\dim\widetilde{R}_e=1$. After this we can invoke Theorem \ref{inv} to obtain the structure of $R=\oplus_{g\in G}\widetilde{R_g}$, as it is claimed in our theorem. To obtain the components of the original grading of $L$, we recall that by Theorem \ref{t00}, $L_{g}=\widetilde{R}_{g}^{(-)}\oplus \widetilde{R}_{gh}^{(+)}\cap L$. Also, $X_t\in L$ for $t\neq e$ and $X_t\in R^{(-)}$ if $\beta(t)=-1$ while $X_t\in R^{(+)}$ if $\beta(t)=1$. Now the computation leading to the explicit form of the components of the original grading becomes obvious.\epf
 

\section{Elementary Gradings}\label{s3}

In this section we continue our study of outer gradings on $L=\rsl(n)$. We use the notation and results of Proposition \ref{np1}. Therefore, $\wh{G}=\langle\vp\rangle\times\Lambda_1$, $\Lambda=\langle\vp^2\rangle\times \Lambda_1$, $G=\langle a\rangle\times K$ so that $\vp(K)=1$, $\vp(a)=\rho$, $\rho$ the $(2m)^\mathrm{th}$ primitive root of $1$, $\Lambda^\perp=\langle h\rangle$ where $h=a^{m}$, $\Lambda_1^\perp=\langle a\rangle$, $\vp(h)=-1$. The Lie automorphism $\vp$ is outer.  By Lemma \ref{lelfin}  the $G/H$ grading of $L$ is inner. We set $\bG=G/H$. Let now the $\bar{G}$-grading on $R=M_n$ be elementary. Suppose it is given by an $n$-tuple
$\tau=(\bg_1,\ldots,\bg_n)$. Then we know from \cite[Section 7]{BShZ} that the action of $\wh{\bG}\cong \Lambda$ is given as follows. If $\lambda\in\Lambda$ and $X\in R$ then
 
\bee{nf11}
\lambda\ast X= \cj{T_\lambda}{X}\mbox{ where }T_\lambda=\mathrm{diag}\,\{\lambda(\bg_1),\ldots,\lambda(\bg_n)\}.
\ene

If we use (\ref{nf5}), then in order for $\vp$ to respect the elementary $\bG$-grading in question, we must have

\bee{nf12}
T_\lambda\Phi T_\lambda=\beta\Phi\mbox{ for all }\lambda\in\Lambda,\; \beta\mbox{ a nonzero scalar, depending on $\lambda$}.
\ene

Now we closely follow the argument in \cite{BShZ} to determine the best possible form of $\Phi$, up to an inner automorphism of $M_n$, making sure that the $\bG$-grading is still elementary. If we identify $M_n$ with $\mathrm{End}\,V$ where $V$ is an $n$-dimensional vector space with basis $\{ e_1,\ldots,e_n\}$ then an elementary grading of $M_n$ corresponding to (\ref{nf11}) is induced by a $\bG$-grading of $V$ where the elements of the basis are homogeneous of degrees forming $\tau$. If we change the places of the elements of the basis then we may assume, without loss of generality, that
$$
\tau=(\underbrace{a_1,\ldots,a_1}_{p_1},\underbrace{a_2,\ldots,a_2}_{p_2},\ldots,\underbrace{a_q,\ldots,a_q}_{p_q})\mbox{ or, shorter, }\tau=(a_1^{(p_1)},\ldots,a_q^{(p_q)})
$$
where $\{ a_1,\ldots,a_q\}=\{\bg_1,\ldots,\bg_n\}$, $a_i\neq a_j$ for $i\neq j$,  $p_1+\ldots+p_q=n$. 

We split $\Phi$ into $q^2$ blocks by drawing horizontal and vertical lines according to the partition $n=p_1+\ldots +p_q$. Let $\Phi_{ij}$ denote the block of dimension $p_i\times p_j$ in the position $(i,j)$. If we use (\ref{nf12}) then we have $\lambda(a_i)\Phi_{ij}\lambda(a_j)=\beta\Phi_{ij}$ for all $1\leq i,j\leq q$. Since $\lambda$ is a  character, if $\Phi_{ij}\neq 0$ then $\lambda(a_ia_j)=\beta$. Now if $\Phi_{ij}\neq 0$ and $\Phi_{i^\prime j}\neq 0$ then $\lambda(a_ia_j)=\lambda(a_{i^\prime}a_j)$ and so $\lambda(a_i)=\lambda(a_{i^\prime})$ for any $\lambda\in\Lambda$. Because $\wh{\bG}\cong \Lambda$ we must have $a_i=a_{i^\prime}$. By our hypothesis then $i=i^\prime$. Recalling that $\Phi$ is non-degenerate we determine that in each row of blocks there is one block $\Phi_{ij}$ different from zero. The same is true for the columns of blocks of $\Phi$. 

Now it follows from the commutativity $a_ia_j=a_ja_i$ that if in the $i^{{\rm th}}$ row we have $\Phi_{ij}\neq 0$ then in the $j^{{\rm th}}$ row we must have $\Phi_{ji}\neq 0$. Our main property used for this and other claims is that the value of $a_ia_j$ is constant for all cases where $\Phi_{ij}\neq 0$. Let us denote by $x_0$ the element of the group $\bG$ to which all such products are equal.  It is obvious now that by rearranging the basis of $V$ and changing notation for the components of $\tau$ we may assume that
$$
\tau=(x_1^{(k_1)},\ldots,x_s^{(k_s)},y_1^{(m_1)},z_1^{(m_1)},\ldots,y_t^{(m_t)},z_t^{(m_t)}).
$$
Also, there exists an element $x_0\in \bG$ such that
$$
x_1^2=\ldots=x_s^2= y_1z_1=\ldots y_tz_t=x_0.
$$
In the same basis we must have
\bee{nf13}
\Phi=\mathrm{diag}\,\left\{\Phi_1,\ldots,\Phi_s,\left[\begin{array}{cc}0&\Psi_1\\ \Psi^\prime_1&0\end{array}\right],\ldots,\left[\begin{array}{cc}0&\Psi_t\\ \Psi^\prime_t&0\end{array}\right]\right\}.
\ene

It should be noted that thanks to the non-degeneracy of $\Phi$ we must have all ingredient matrices square and $\Psi_i$ of the same order as $\Psi^\prime_i$.

Now we have to recall $\vp^2\in\Lambda$. Let us set $\lambda_0=\vp^2$. Then $\vp^2\ast X=\cj{T_{\lambda_0}}{X}$, for any $X\in L$. Then we should have

$$
\vp^2\ast X=\cj{\Phi}{\tr{(\cj{\Phi}{\tr{X}})}}=\cj{(\tr{\iv{\Phi}}\Phi)}{X}.
$$

It follows that 
\bee{nf14}
\tr{\iv{\Phi}}\Phi=\alpha T_{\lambda_0},\mbox{ that is, }\Phi=\alpha\tr{\Phi}T_{\lambda_0}.
\ene

Considering the explicit form (\ref{nf13}) for $\Phi$, we find that the following are true:

\bee{nf15}
\Phi_i=\alpha\lambda_0(x_i)\tr{\Phi}_i\mbox{ for all }1\leq i\leq s.
\ene
and
\bee{nf16}
\left[\begin{array}{cc}
        0  & \Psi_j   \\ 
   \Psi_j^\prime       &  0  \\ 
\end{array}\right]
=\left[\begin{array}{cc}
       0   & \tr{\Psi_j^\prime}   \\ 
       \tr{\Psi_j}   &  0  \\ 
\end{array}\right]\left[\begin{array}{cc}
         \alpha\lambda_0(y_j) & 0   \\ 
      0    &   \alpha\lambda_0(z_j)  \\ 
\end{array}\right].
\ene

It follows from (\ref{nf15}) that $\alpha^2\lambda_0(x_i)^2=\alpha^2\lambda_0(x_0)=1$. Therefore each $\Phi_i$ is either symmetric or skew-symmetric. We also have
\bee{epsi}
\Psi_j=\alpha\lambda_0(z_j)\tr{\Psi_j^\prime}=\alpha^2\lambda_0(z_j)\lambda_0(y_j)\Psi_j.
\ene
So again $\alpha^2\lambda_0(z_j)\lambda_0(y_j)=\alpha^2\lambda_0(x_0)=1$. 

Now let us conjugate our matrix algebra by a nondegenerate matrix of a shape similar to $\Phi$:

$$
P=\mathrm{diag}\,\left\{A_1,\ldots,A_s,\left[\begin{array}{cc}B_1&0\\ 0&C_1\end{array}\right],\ldots,\left[\begin{array}{cc}B_t&0\\ 0&C_t\end{array}\right]\right\}
$$

In this case the matrices $T_\lambda$ will not change while the blocks of $\Phi$ will change as follows. By Lemma \ref{nl1}  each $\Phi_i$ will be replaced by $\tr{A_i}\Phi_iA_i$. This allows us to assume, without loss of generality, that, depending on whether $\alpha\lambda_0(x_i))$ is $1$ or $-1$, each $\Phi_i$ is either an identity matrix $I_{k_i}$ of appropriate order $k_i$  or  a matrix $S_{l_p}$, $S_{l_p}=\left[\begin{array}{cc}
         0 &  I_{l_p} \\ 
         -I_{l_p} &  0   
\end{array}\right]$ of even order $k_p=2l_p$.

Now the block $\left[\begin{array}{cc}0&\Psi_j\\ \Psi_j^\prime&0\end{array}\right]$ will be replaced by $\left[\begin{array}{cc}0&\tr{B_j}\Psi_jC_j\\ \tr{(\tr{B_j}\Psi_j^\prime C_j)}&0\end{array}\right]$. Considering (\ref{epsi}), this allows us to assume, without loss of generality, that the whole block in question can be replaced by  $\left[\begin{array}{cc}0&I_{m_u}\\ \alpha\lambda_0(y_u)I_{m_u}&0\end{array}\right]$ where $I_{m_u}$ is the identity matrix of an appropriate order $m_u$.

In the new basis $\Phi$ will look like the following:

\begin{eqnarray}
\Phi&=&\mathrm{diag}\,\left\{I_{k_1},\ldots, I_{k_r},\left[\begin{array}{cc}0&I_{l_{r+1}}\\ -I_{l_{r+1}}&0\end{array}\right],\ldots, \left[\begin{array}{cc}0&I_{l_{s}}\\ -I_{l_s}&0\end{array}\right]\right.,\nonumber\\ &&\left. \left[\begin{array}{cc}0&I_{m_1}\\ \alpha\lambda_0(y_1)I_{m_1}&0\end{array}\right],\ldots,\left[\begin{array}{cc}0&I_{m_t}\\ \alpha\lambda_0(y_t)I_{m_t}&0\end{array}\right]\right\}, \label{nf1705}
\end{eqnarray}

At this point we are ready to construct an inner automorphism $\psi$ of $L$ such that $\psi$ commutes with the action of the whole of $\wh{G}$ and also the action of $\vp^2=\lambda_0$ coincides with $\psi^2$, which will allow us to apply Theorem \ref{t00}.

We will look for $\psi$ in the form of an inner automorphism given by a diagonal matrix $T_\psi$ with respect to a basis of $V$ which results after all the above transformations. Notice that thanks to (\ref{nf14}) and (\ref{nf1705}), in which we set $\gamma_u=\alpha\lambda_0(y_u)$ for $u=1,\ldots,t$, we can write

\begin{eqnarray}
T_{\lambda_0}&=&\alpha^{-1}\,{\tr{\Phi}}^{-1}\Phi=\alpha^{-1}\mathrm{diag}\left\{I_{k_1},\ldots,I_{k_s},-I_{l_{r+1}},-I_{l_{r+1}}\ldots,-I_{l_s},-I_{l_s},\right.\nonumber\\ &&\left.\gamma_1 I_{m_1},\gamma_1^{-1}I_{m_1},\ldots, \gamma_tI_{m_t},\gamma_t^{-1}I_{m_t}\right\}.\label{nf1702}
\end{eqnarray}

This suggests that the matrix we want to find has the form of

\begin{eqnarray}
T_{\psi}&=&\mathrm{diag}\,\left\{\ve I_{k_1},\ldots,\ve I_{k_s},\pi I_{l_{r+1}},\rho I_{l_{r+1}}\ldots,\pi I_{l_s},\rho I_{l_s},\right.\nonumber\\ &&\left.\mu_1 I_{m_1},\nu_1 I_{m_1},\ldots, \mu_t I_{m_t},\nu_t I_{m_t}\right\}.\label{nf1703}
\end{eqnarray}

Now $T_\psi^2=\xi T_{\lambda_0}$ for some scalar $\xi$. This gives us the following relations

\begin{eqnarray}
\ve^2&=&\xi\alpha^{-1},\nonumber\\
\pi^2&=&-\xi\alpha^{-1},\,\rho^2=-\xi\alpha^{-1},\nonumber\\
\mu_u^2&=&\xi\gamma_u,\,\nu_u^2=\xi\gamma_u^{-1},\, u=1,\ldots,t.\label{nf1707}
\end{eqnarray}

Now in order for $\psi$ to commute with $\vp$ we must have an equation of the form (\ref{nf12}) satisfied for $T_\psi$, that is,

\bee{nf1704}
T_\psi\Phi T_\psi=\delta\Phi,
\ene
for an appropriate parameter $\delta$. If we use (\ref{nf1704}) then, considering all nonzero entries of $\Phi$,
we arrive at the following relations

\begin{eqnarray}
\ve^2&=&\delta,\nonumber\\
\pi\rho&=&\delta,\nonumber\\
\mu_u\nu_u&=&\delta.\label{nf1706}
\end{eqnarray}

In the case where at least one of the diagonal blocks of $\Phi$ is nonzero, we can modify our $n$-tuple $\tau$ by dividing all entries by $x_1$. Then we will have $x_0=\be$ and $\alpha^2=1$. Resolving (\ref{nf1707}) and (\ref{nf1706}) in this case we obtain $\xi\alpha^{-1}=\delta$ and then $\pi\rho=\xi\alpha^{-1}$ and $\mu_u\nu_u=\xi\alpha^{-1}$, for all $u$. It follows that, in this case, one can choose $\ve$ any square root of $\xi\alpha^{-1}$, $\pi$ any square root of $-\xi\alpha^{-1}$  while $\rho=-\pi$. Then $\pi\rho=\xi\alpha^{-1}$, as required. We can also take as $\mu_u$ any square root of $\xi\gamma_u$  and $\nu_u=\frac{\xi\alpha}{\mu_u}$. Then $\nu_u^2=\xi\gamma_u^{-1}$, as needed. At the same time, $\mu_u\nu_u=\xi\alpha=\xi\alpha^{-1}$ because we have assumed $\alpha^2=1$.

In the case where we have no nonzero diagonal blocks in $\Phi$ the only equations we have to resolve are
$$
\mu_u^2=\xi\gamma_u,\,\nu_u^2=\xi\gamma_u^{-1}\mbox{ and }\mu_u\nu_u=\delta.
$$
In this case we take $\mu_u$ any square root of $\xi\gamma_u$ and $\nu_u=\frac{\xi\alpha}{\mu_u}$. Then the value of $\mu_u\nu_u$ is a constant $\xi\alpha$ so that we can set $\delta$ equal to this value.

Now $T_\psi$ satisfies (\ref{nf1707}) hence $\psi$ commutes with $\vp$. Since $T_{\psi}$ is diagonal the conjugation by $T_{\psi}$ commutes with the conjugation by all $T_\lambda$, $\lambda\in\Lambda$.

Thus the existence of $\psi$ with the properties desired has been proved. Once we found $\psi$ we can use Theorem \ref{t00} and make the following conclusion.

\begin{theorem}\label{telem} Let $L=\rsl(n)$ be given an outer grading by a finite abelian group $G$ such that the respective $\bG$-grading is elementary. Then there is an involution $\ast$ on $R=M_n$, an element $h$ of order  $2$ in $G$, and an elementary involution $G$-grading  $R=\sum_{g\in G}\widetilde{R}_g$ such that
$$
L_g=\left\{\begin{array}{ll}\widetilde{R}_g^{(-)}\oplus \widetilde{R}_{gh}^{(+)}&\mbox{ if }g\neq h\\
\widetilde{R}_h^{(-)}\oplus (\widetilde{R}_{e}^{(+)}\cap L)&\mbox{ otherwise}.\end{array}\right.
$$
Here $\widetilde{R}_g^{(\pm)}$ is the set of symmetric (skew-symmetric) elements in $\widetilde{R}_g$ with respect to the involution $\ast$.
\end{theorem}

\pp All claims have been proved except that the grading $R=\sum_{g\in G}\widetilde{R}_g$ is elementary. This, however, easily follows because the action of $\psi$ is given as the conjugation by a diagonal matrix. So every matrix unit $E_{ij}$ is an eigenvector of $\psi$. Since the same is true for the action of $\Lambda_1$, every matrix unit $E_{ij}$ is graded. It is well-known in this case  \cite{ZB} that the grading must be elementary.\epf

\section{Mixed Gradings}

Our approach to handling the outer gradings on $L=\rsl(n)$, is to apply Theorem \ref{t00}. Therefore, given an outer automorphism of $L=\rsl(n)$, we have to find an inner automorphism $\psi$ which commutes with the action of $\wg$ and such that the action of $\vp^2$ coincides with $\psi^2$. 

At this point we have $R=R^{(0)}\otimes R^{(1)}\otimes\cd\otimes R^{(k)}$ where the $\bG$-grading on $R^{(0)}$ is elementary and that on $R^{(1)}\otimes\cd\otimes R^{(k)}$ is fine. Recall also that $G$-grading on Lie algebra $L$ induces a $G$-grading on $R$ as on
    a Lie algebra. The subspaces $R^{(0)}$ and $R^{(1)}\otimes\cd\otimes R^{(k)}$ are also $G$-graded and the $G$-grading on them has  been described in Sections \ref{s2} and \ref{s3}. Let us choose $G$-graded bases in these two subspaces, such that every element is either traceless or is the identity matrix. Then we have a $G$-graded basis on $R$ consisting of $I\otimes I$ and some traceless matrices of the form $u\otimes v$ where at least one of $u,v$  has trace zero. Then these latter matrices  will form a $G$-graded basis of $L=\rsl(n)$. To prove this claim we have to apply the generators of $\wh{G}$ to these matrices. Assume $g_1=\deg_G u$, $g_2=\deg_G v$. If none of $u$, $v$ is $I$, and $\lambda\in \Lambda$ we will have 
$$
\lambda\ast (u\otimes v)=(\lambda\ast u)\otimes (\lambda\ast v)=\lambda(g_1)\lambda(g_2)(u\otimes v)=\lambda(g_1g_2)(u\otimes v)=\lambda(g_1g_2h)(u\otimes v).
$$
If we apply $\vp$ then using (\ref{nf6}) we will get
$$
\vp\ast (u\otimes v)=-(\vp\ast u)\otimes (\vp\ast v)=-\vp(g_1)\vp(g_2)(u\otimes v)=\vp(g_1g_2h)(u\otimes v).
$$
It follows that $\deg(u\otimes v)=g_1g_2h$. As for the elements of the form $u\otimes I$ and $I\otimes v$ then they retain their degrees as the elements of $R^{(0)}$ or $R^{(1)}\otimes\cd\otimes R^{(k)}$, that is, $\deg(u\otimes I)=g_1$, $\deg(I\otimes v)=g_2$.

Using this notation, we define a mapping $\psi:R\rightarrow R$ by the formula
\bee{mg1} \psi\ast (u\otimes v)= (\psi\ast u)\otimes v.
\ene
Here $\psi\ast u$ is defined as an inner automorphism, the result of our argument in Section \ref{s3} leading to Theorem \ref{telem}. For that $\psi$ we had $\psi^2=\vp^2$ on $R^{(0)}$.
We also remember that according to Corollary \ref{fin2}, $\vp^2=\mathrm{id}$ when restricted to $R^{(1)}\otimes\cd\otimes R^{(k)}$. So it is immediate $\psi^2=\vp^2$ for the mapping defined by (\ref{mg1}). Clearly, $\psi$ is inner, given by the matrix $T_\psi\otimes I$ where $T_\psi$ has been found in Section \ref{s3}. It is obvious that $\psi$ commutes with $\wh{G}$. Thus our $G$-grading of $L=\rsl(n)$ can be recovered by Theorem \ref{tomega} from a $G$-grading, which respects an involution of $R=M_n$. All such gradings have been completely described in Theorem \ref{inv}.

Our final results will then look as follows. 

\begin{theorem}\label{tmg} Let $F$ be an algebraically closed field of characteristic zero. Any grading $L=\bigoplus_{g\in G}L_g$ of $L=\rsl(n)$ by a finite abelian group $G$ on $R=M_n$ is conjugate by an inner automorphism of $R$ to one of the following types.
\begin{enumerate}
\item[{\rm (I)}] The restriction to $L$ of any associative $G$-grading of $M_n$;
\item[{\rm (II)}] Given any involution $G$-grading $R=\bigoplus_{g\in G}\widetilde{R}_g$ and an element $h$ of order $2$ in $G$, the grading defined by
$$
L_g=\left\{\begin{array}{ll}\widetilde{R}_g^{(-)}\oplus \widetilde{R}_{gh}^{(+)}&\mbox{ if }g\neq h\\
\widetilde{R}_h^{(-)}\oplus (\widetilde{R}_{e}^{(+)}\cap L)&\mbox{ otherwise}.\end{array}\right.
$$
\end{enumerate}
\end{theorem} 

\bigskip

If we want a more explicit form, we have to notice that in both ``inner'' and ``outer'' types of the gradings we have to start with a $G$-grading $R=A\otimes B$ where $A\cong M_p$ is a $G$-graded subalgebra with an elementary $G$-grading and $B\cong M_q$ a $G$-graded subalgebra with a fine grading. There is a subgroup $T\subset G$, $T\cong \mathbb{Z}_{n_1}^2\times\ld\times\mathbb{Z}_{n_k}^2$ with $T\cap \su{M_p}=\{ e\}$, which supports $M_q$. Thus  a basis of $B$ can be chosen in the form of $\{ X_t\,|\,t\in T\}$ where each $X_t$ is the Kronecker product $X_{t_1}\otimes\ld \otimes X_{t_k}$ where $t_s=a_s^{i_s}b_s^{j_s}\in \mathbb{Z}_{n_s}^2$ and $X_{t_s}=X_{a_s}^{i_s}X_{b_s}^{j_s}$, as in (\ref{e00}), $s=1,\ld,k$. Also, there is a $p$-tuple $\tau=(g_1,\ld,g_p)$ of elements of $G$, which defines the elementary grading of $A$ by $\deg E_{ij}=g^{-1}_ig_j$. 

In the case of Type I gradings, to obtain any grading of $L$ we have to choose a basis of $R$ in the form $\{ E_{ij}\otimes X_t\,|\, 1\leq i,j\leq p, t\in T\}$ and set 
\beq{eTI}
L_g&=&\mathrm{Span}\left\{\left\{ E_{ij}\otimes X_t\,|\, g^{-1}_ig_jt=g, \, 1\leq i\neq j\leq n\right\}\right.\\
&&\left.\bigcup \left\{(E_{11}-E_{ii})\otimes X_g\,|\, 1<i\leq n\right\}\right\}.\nonumber
\eqe 

In the case of Type II gradings, the $p$-tuple $\tau$ has to be chosen as prescribed in Lemmas \ref{L8'} and \ref{L8}. The same lemmas define an involution $\ast$ on $A$. Also, $n_1=\ld=n_k=2$,  and an involution $\ast$ is defined on $B$ by $X_t^\ast=-(\mathrm{sgn}\,t)X_t$ where $\mathrm{sgn}\,t=\pm 1$, as defined in Theorem \ref{inv}. Now an involution $\ast$ is defined on $R$ by $(Y\otimes X_t)^\ast=Y^\ast\otimes X_t^\ast$ so that $Y\otimes X_t$ is symmetric if either $Y$ is symmetric and $\mathrm{sgn}\,t=1$ or $Y$ is skew-symmetric and  $\mathrm{sgn}\,t=-1$. Similarly for skew-symmetric elements. It is obvious from the disjoint property for the supports of $A$ and $B$ that such symmetric (skew-symmetric) elements span $R^{(\pm)}$. Now we can define

\beq{eTIIg}
L_g&=&\mathrm{Span}\,\left\{\left\{ Y\otimes X_t\,|\, (\deg Y)t=g, Y\otimes X_t\mbox{ skew-symmetric}\right\}\right.\\&&\left.\bigcup\left\{ Y\otimes X_t\,|\, (\deg Y)t=gh, Y\otimes X_t\mbox{ symmetric}\right.\}\right\}\nonumber
\eqe
if $g\neq h$ and
\beq{eTIIh}
L_h&=&\mathrm{Span}\,\left\{\left\{ Y\otimes X_t\,|\, (\deg Y)t=h, Y\otimes X_t\mbox{ skew-symmetric}\right\}\right.\\ &&\bigcup\left\{ Y\otimes X_t\,|\, (\deg Y)t=e, Y\otimes X_t\mbox{ symmetric,}\right.\nonumber\\&&\left.\left.\mbox{ and }\mathrm{Tr}\, Y=0\mbox{ if }t=e\right\}\right\}.\nonumber
\eqe

\end{document}